%% file: CharactSome.tex
\documentclass[12pt,reqno]{amsart}
\usepackage{amsthm}
\usepackage{amsmath}
\usepackage{version}
\usepackage{graphicx}
\usepackage{enumerate}
\usepackage[usenames,dvipsnames] {color}
\numberwithin{equation}{section}

\def\d{\mathbb{D}}
\def\c{\mathbb{C}}
\def\t{\mathbb{T}}

\def\cald{\mathcal{D}}

\newcommand\calc{\mathcal{C}}

\def\car#1{|#1|_{\rm car}}
\def\kob#1{|#1|_{\rm kob}}

\def\be{\begin{equation}}
\def\ee{\end{equation}}

\def\s0{s_0}
\def\p0{p_0}

\newcommand\nd{nondegenerate }

\DeclareMathOperator{\Car}{{\mathrm Car}}
\DeclareMathOperator{\Kob}{{\mathrm Kob}}

 \newtheorem{theorem}{Theorem}[section]
 \newtheorem{corollary}[theorem]{Corollary}
 \newtheorem{lemma}[theorem]{Lemma}

\newtheorem{example}[theorem]{Example}
\newtheorem{fact*}{Fact}

\DeclareMathOperator\re{Re}

\DeclareMathOperator\hol{Hol}

\newcommand\half{{\tfrac 12}}

\newcommand\idd{\mathrm{id}_\mathbb{D}}

\newcommand{\T}{\mathbb{T}}

\newcommand{\D}{\mathbb{D}}
\newcommand{\C}{\mathbb{C}}

\newcommand{\set}[1]{\left\{#1\right\}}
\newcommand{\norm}[1]{\left\Vert#1\right\Vert}

\newcommand{\ran}[1]{\operatorname{ran}#1}

\newcommand{\inv}{^{-1}}

\newcommand{\ph}{\varphi}
\renewcommand\phi{\varphi}
\newcommand\al{\alpha}
\newcommand\ga{\gamma}

\newcommand\de{\delta}

\newcommand\la{\lambda}

\newcommand\beq{\begin{equation}}

\newcommand\eeq{\end{equation}}
\newcommand\df{\stackrel{\rm def}{=}}

\newcommand\bbm{\begin{bmatrix}}
\newcommand\ebm{\end{bmatrix}}
\newcommand\bpm{\begin{pmatrix}}
\newcommand\epm{\end{pmatrix}}
\numberwithin{equation}{section}

\newtheorem{lem}[theorem]{Lemma}
\newtheorem{prop}[theorem]{Proposition}
\newtheorem{thm}[theorem]{Theorem}

\newtheorem{defin}[theorem]{Definition}

\begin{document}
\title[Characterization of some domains]{Characterizations of some domains via Carath\'eodory extremals}
\author{Jim Agler}
\address{Department of Mathematics, University of California at San Diego, CA \textup{92103}, USA}
\thanks{Partially supported by National Science Foundation Grants
DMS 1361720 and 1665260, a Newcastle URC Visiting Professorship,  London Mathematical
Society Grant 41730  and the Engineering and Physical Sciences Research Council grant EP/N03242X/1 }

\author{Zinaida Lykova}
\address{School of Mathematics, Statistics and Physics, Newcastle University, Newcastle upon Tyne
 NE\textup{1} \textup{7}RU, U.K.}
\email{Zinaida.Lykova@ncl.ac.uk}

\author{N. J.  Young}
\address{School of Mathematics, Statistics and Physics, Newcastle University, Newcastle upon Tyne NE1 7RU, U.K.
{\em and} School of Mathematics, Leeds University,  Leeds LS2 9JT, U.K.}
\email{Nicholas.Young@ncl.ac.uk}
\date{21st February 2018, revised 15th July 2018}
\begin{abstract} 
In this paper we characterize the unit disc, the bidisc and the symmetrized bidisc
\[
G =\{(z+w,zw):|z|<1,\ |w|<1\}
\]
in terms of the possession of small classes of analytic maps into the unit disc that suffice to solve all Carath\'eodory extremal problems in the domain.
\end{abstract} 

\maketitle

\input intro

\input CandK

\input charactbidisc

\input charGcara

\input uniqueness

\input concluding

\bibliography{references}

\end{document}

%% file: intro.tex
\section*{Introduction}
The  Carath\'eodory extremal functions on a domain $\Omega\subset\c^d$ constitute a special class of analytic maps from $\Omega$ into the unit disc $\d$, a class which contains a lot of information about $\Omega$.  
By a {\em domain} we mean a connected open set in $\c^d$ for some integer $d\geq 1$.   

In brief, the Carath\'eodory extremal functions on $\Omega$ are those which, for some pair $z,w$ of distinct points of $\Omega$, maximise over all analytic maps $F$ from $\Omega$ to the open unit disc $\d$ the Poincar\'e distance from $F(z)$ to $F(w)$.  A function $F$ for which the maximum is attained is said to {\em solve the Carath\'eodory extremal problem} for the pair $z,w$.

Even for such a simple domain as the bidisc $\d^2$ the class of all Carath\'eodory extremal functions is large, and there are few domains for which an explicit description of all Carath\'eodory extremal functions can be given.  However, it can happen that there is a relatively small set $\calc$ of Carath\'eodory extremal functions on a domain $\Omega$ which is {\em universal for the Carath\'eodory extremal problem}, in the sense that, for every choice of distinct points $z,w$ in $\Omega$, there exists a function $F\in\calc$ that solves the Carath\'eodory extremal problem for  $z,w$.
See Definition \ref{predef10} for a more formal statement.

For example, if $\Omega=\d$ then we may take $\calc$ to comprise only the identity map on $\d$, while if $\Omega=\d^2$ then the set $\calc$ comprising the two co-ordinate functions has the universal property.

 In this paper we pose the question {\em to what extent does knowledge of a universal set for the Carath\'eodory extremal problem on $\Omega$ determine $\Omega$ up to isomorphism?} 

One has to make some assumption to rule out cases such as $\Omega=\c^d$, in which the only analytic functions on $\Omega$ are constant functions.
We therefore restrict attention to Lempert domains (see Definition \ref{predef50}),  for which the theory of hyperbolic complex spaces (in the sense of \cite{kob98})  is suitably rich.

We give a positive answer to the question for three domains, namely the disc $\d$, the bidisc $\d^2$ and the {\em symmetrized bidisc}
\beq\label{defG}
G\df\{(z+w,zw):|z|<1,\ |w|<1\} \, \subset \c^2.
\eeq
For each of these domains there is a small class of functions which is universal for the  Carath\'eodory extremal problem on the domain, and moreover an appropriate structure of this small class does determine the isomorphism class of the domain among Lempert domains (Theorems \ref{chardisc}, \ref{thm10} and \ref{carthm10}).  In the first two instances, the structure is just the cardinality, while in the case of $G$, it is a special linear fractional parametrization by the unit circle $\t$.

The domain $G$ is of interest in connection with the theory of invariant distances \cite{jp}.  $G$ has an extensive literature, including 
\cite{ay2004,cos04,ez05,jp04,pz05,kos,tryb,aly2016,bhatta,sarkar}. 

In Section \ref{CandK} we describe the Carath\'eodory and Kobayashi extremal problems, together with basic facts about complex geodesics.
In Section \ref{charactbidisc} we characterize the disc and bidisc up to isomorphism among Lempert domains by the property that they admit minimal universal sets for the Carath\'eodory problem comprising one and two functions respectively.

In the symmetrized bidisc there is a one-parameter family of functions which is universal for the Carath\'eodory problem on $G$ (\cite[Theorem 1.1 and Corollary 4.3]{ay2004}).  This family comprises the rational functions
\be\label{pre20_0}
\Phi_\omega(s^1,s^2) = \frac{2\omega s^2 -s^1}{2-\omega s^1},
\ee
where $\omega \in \t$. 

In Section \ref{charactcara} 
 we prove a converse to this statement:  a domain $\Omega$ in $\c^2$ is biholomorphic to $G$ if and only if $\Omega$ is a Lempert domain and $\Omega$ has a universal set for the Carath\'eodory problem that admits a  linear fractional parametrization by $\t$ akin to equation \eqref{pre20_0} (Theorem \ref{carthm10}). 

In Section \ref{uniqueness} we show that, in the case of the three domains $\d$, $\d^2$ and $G$, the minimal universal set for the Carath\'eodory problem is {\em unique} up to a natural notion of equivalence.

We are grateful to a referee for some useful comments.

%% file: CandK.tex
\section{The Carath\'eodory and Kobayashi problems}\label{CandK}
In this section we describe our terminology for the Carath\'eodory and Kobayashi extremal problems.

If $U \subseteq \c^{n_1}$ and $V \subseteq \c^{n_2}$ are two open sets we denote by $V(U)$ the set of holomorphic mappings from $U$ into $V$.

If $U$ is an open set in $\c^n$, then by a \emph{datum in} $U$
 we mean an ordered  pair $\delta$ where either $\delta$ is \emph{discrete}, that is, has the form
\[
\delta =(s_1,s_2)
\]
where $s_1,s_2 \in U$, or $\delta$ is \emph{infinitesimal}, that is, has the form
\[
\delta = (s,v)
\]
where $s \in U$ and $v\in\c^n$. 

If $\delta$ is a datum, we say that $\delta$ is \emph{degenerate} if either $\delta$ is discrete and $s_1=s_2$ or $\delta$ is infinitesimal and $v=0$. Otherwise, we say that $\delta$ is \emph{nondegenerate}.

An infinitesimal datum in $U$ is the same thing as a point of the complex tangent bundle $TU$ of $U$.

For $F\in \Omega(U)$, $s\in U$, and $v \in \c^{n_1}$, the directional derivative $D_v F(s) \in \c^{n_2}$ is defined by
\index{$D_vF$}
\[
D_v F(s) = \lim_{z \to 0} \frac{F(s+zv) - F(s)}{z}.
\]

 If $U$ and $\Omega$ are domains, $F\in \Omega(U)$, and $\delta$ is a datum in $U$, we define a datum $F(\delta)$ in $\Omega$ by
 \[
 F(\delta) = (F(s_1),F(s_2))
 \]
 when $\delta$ is discrete and by
 \[
 F(\delta) =(F(s),D_v F (s))
 \]
 when $\de$ is infinitesimal. 

 For any datum $\delta$ in $\d$, we define $|\delta|$ to be the Poincar\'e distance or metric at $\de$ in the discrete or infinitesimal case respectably, that is
\[
|\delta| = \tanh\inv \left|\frac{z_1 -z_2}{1-\bar{z_2}z_1}\right|
\]
when $\delta=(z_1,z_2)$ is discrete\footnote{In \cite[Chapter 3]{aly2016} we used a different notation in that we omitted $\tanh\inv$; this makes no essential difference.},
and by
\[
|\delta|=\frac{|v|}{1-|z|^2}
\]
when $\delta = (z,v)$ is infinitesimal.
\index{$|\delta|$}
\\ 

{\bf  The Carath\'eodory extremal problem.}
\index{Carath\'eodory extremal problem}
{\em For a domain  $U$ in $\c^n$ and a nondegenerate datum $\delta$  in $U$, compute the quantity $\car{\delta}$
defined by}
\[
\car{\delta}=\sup_{F\in \d(U)} |F(\delta)|.
\]
We shall refer to this problem as $\Car\delta$ 
and will say that \emph{$C$ solves} $\Car \delta$ if $C\in \d(U)$ and
\[
\car{\delta} = |C(\delta)|.
\]
When it is necessary to specify the domain in question we shall write $\car{\de}^U$.

It is easy to see with the aid of Montel's theorem that, for every \nd datum $\de$ in $U$, there does exist $C\in\d(U)$ which solves $\Car\de$.  Such a $C$ is called a {\em Carath\'eodory extremal function } for $\de$.

\begin{defin}\label{predef10} 
For  a domain $U$ in $\c^n$, we say that a set $\mathcal{C} \subseteq \d(U)$ is a \emph{universal set for the Carath\'eodory extremal problem on $U$}
 if, for every nondegenerate datum $\delta$  in $U$, there exists a function $C \in \mathcal{C}$ such that $C$ solves $\Car \delta$.
If, furthermore, no proper subset of $\mathcal C$ is universal for the Carath\'eodory extremal problem on $U$, then we say that $\mathcal C$ is a {\em minimal} universal set. 
\end{defin}

\begin{example}\label{ex1.2} {\rm For  many classical domains $\Omega$ there are small  sets which are universal  for the Carath\'eodory extremal problem on $\Omega$. 

(i) If $\Omega = \D^d$ then the set of the $d$ co-ordinate functions is a minimal universal set for the Carath\'eodory extremal problem on $\D^d$, as follows from Schwarz' Lemma. 

(ii) For the open Euclidean unit ball $\mathbb{B}_d$ in $\C^d$, there is a universal set for the Carath\'eodory extremal problem consisting of compositions of the projections onto the planes through the center with automorphisms of $\mathbb{B}_d$.

(iii) If $\Omega = G$, the symmetrized bidisc, there is a $1$-parameter set $\{ \Phi_{\omega}: \omega \in \T \}$  (see equation \eqref{pre20} below) that constitutes  a minimal universal set for the Carath\'eodory extremal problem on  $G$ \cite{ay2004}. 
}
\end{example}

\begin{defin}\label{predef20}
We say that a domain $U$ in $\c^n$ is \emph{weakly hyperbolic}
\index{weakly hyperbolic}
 if $\car{\delta} >0$ for every nondegenerate datum $\delta$ in $U$. Equivalently, for every nondegenerate datum $\delta$ in $U$, there exists a bounded holomorphic function $F$ on $U$ such that $F(\delta)$ is a nondegenerate datum in $\c$.
\end{defin}
\begin{lem}\label{prelem10}
Let $\Omega$ be a weakly hyperbolic domain in $\c^n$ and assume that $\mathcal{C}$ is a universal set for the Carath\'eodory extremal problem on $\Omega$. If $\delta$ is a nondegenerate datum in $\Omega$, then there exists $C\in\mathcal{C}$ such that $C(\delta)$ is a nondegenerate datum in $\d$.
\end{lem}
\begin{proof}
If $C(\delta)$ is degenerate for all $C\in\mathcal{C}$, then
\begin{align*}
\car{\delta} =\sup_{C \in \mathcal{C}}  |C(\delta)| =0.
\end{align*}
But as $\Omega$ is assumed to be weakly hyperbolic and $\delta$ is assumed to be nondegenerate,
\[
\car{\delta} >0.
\]
\end{proof}

\noindent{\bf  The Kobayashi extremal problem.} 
\index{Kobayashi extremal problem}
{\em For  a domain $U$ in $\c^n$ and a nondegenerate datum $\delta$ in $U$, compute the quantity $\kob{\delta}$
 defined by}
\beq\label{defkob}
\kob{\delta}=\inf_{\substack{f\in U(\d)\\ f(\zeta) =\delta}} |\zeta|.
\eeq
Again, where appropriate we shall indicate the domain by a superscript.
We shall refer to this problem as $\Kob \delta$ 
and will say that \emph{$k$ solves} $\Kob\delta$
 if $k\in U(\d)$ and there exists a datum $\zeta$ in $\d$ such that $k(\zeta) = \delta$ and
\[
\kob{\delta} = |\zeta|.
\]

Note that the infimum in the definition \eqref{defkob} of $\kob{\de}$ is attained if $U$ is a taut domain, where $U$ is said to be {\em taut} if $\hol(\d,U)$ is a normal family \cite[Section 3.2]{jp}.   In particular $\kob{\de}$ is attained when $U=G$  \cite[Section 3.2]{jp}.  Any function which solves $\Kob\de$ is called a {\em Kobayashi extremal function for $\de$}.

Let $U$ be a domain in $\c^n$ and $\delta$ a nondegenerate datum in $U$. The solutions to $\Car \delta$ and $\Kob\delta$ are never unique, for if $m$ is a M\"obius transformation of $\d$, then $m\circ C$ solves $\Car\delta$ whenever $C$ solves $\Car\delta$ and $f \circ m$ solves $\Kob\delta$ whenever $f $ solves $\Kob\delta$. This suggests the following definition.  
\begin{defin}\label{predef30}
Let $U$ be a domain in $\c^n$ and let $\delta$ be a nondegenerate datum in $U$. We say that \emph{the solution to $\Car\delta$ is essentially unique}, if whenever $F_1$ and $F_2$ solve $\Car\delta$ there exists a M\"obius transformation $m$ of $\d$ such that $F_2=m \circ F_1$.
We say that \emph{the solution to $\Kob\delta$ is essentially unique} if 
 the infimum in equation \eqref{defkob} is attained and, 
 whenever $f_1$ and $f_2$ solve $\Car\delta$ there exists a M\"obius transformation $m$ of $\d$ such that $f_2=f_1 \circ m$.
\end{defin}

\subsection{Complex geodesics}\label{complexgeos}
One of the most striking results in the theory of hyperbolic complex spaces is Lempert's theorem \cite{lem81}
 to the effect that $\car{\cdot}=\kob{\cdot}$ in bounded convex domains.  A consequence is that,
if $\delta$ is a datum in a bounded convex domain $U \subseteq \c^n$, then there exists a solution $C$ to $\Car\delta$ and a solution $k$ to $\Kob\delta$ such that
\be\label{pre10}
C \circ k =\idd.
\ee
In the event that $k:\d \to U$ and there exists $C:U \to \d$ such that equation \eqref{pre10} holds, then necessarily $C$ solves  $\Car\delta$  and $k$ solves $\Kob\delta$. 
In this case $\ran k$ is a complex geodesic.

\begin{defin}
Let $U$ be a domain in  $\c^n$ and let 
$\mathcal{D} \subset U$. We say that $\mathcal{D}$ is a {\em complex geodesic} in $U$ if there exists a function 
$k \in U(\d)$ and a function $C \in \d(U)$ such that 
$C \circ k =\idd$ and  $\mathcal{D} =k(\d)$. 
\end{defin}

 This terminology is suggested by the fact that if $k:\d \to U$, then $\ran k$ is a totally geodesic one-complex-dimensional manifold properly embedded in $U$ if and only if there exists $C:U\to \d$ such that equation \eqref{pre10} holds.

The following definition provides language to  describe in a concise way the relationship between datums in $U$ and complex geodesics in $U$.
 \begin{defin}\label{predef40}
 If $V\subseteq U$ and $\delta$ is a datum in $U$ we say that $\delta$ \emph{contacts} $V$ 
\index{contacts}
if either $\delta=(s_1,s_2)$ is discrete and  $s_1\in V$ and $s_2 \in V$, or $\delta =(s,v)$ is infinitesimal and there exist two sequences of points $\{s_n\}$  and $\{t_n\}$ in $V$ such that $s_n \neq t_n$ for all $n$, $s_n \to s$ and
 \[
 \frac{t_n-s_n}{\norm{t_n-s_n}} \to v_0,
 \]
where $v \sim v_0$.
 \end{defin}
 Note that if $U$ and $\Omega$ are domains, $V$ is a subset of $U$, $\delta$ is a datum in $U$ that contacts $V$ and $F_1$ and $F_2$ are any two holomorphic mappings from $U$ to $\Omega$ satisfying $F_1|V=F_2|V$ then $F_1(\delta) =F_2(\delta)$ \cite[Remark 4.3]{aly2016}.

The following proposition relates the concept of contact to the notion of complex geodesics \cite[Proposition 4.4]{aly2016}.
\begin{prop}\label{preprop10}
Let $U$ be a domain in $\c^n$ and let $k \in U(\d)$ be such that 
 $\mathcal{D} =k(\d)$ is a complex geodesic in $U$. 
\begin{enumerate}[\rm (1)]
\item If $\zeta$ is a datum in $\d$, then $k(\zeta)$ contacts $\cald$. 
\item If $\delta$ is a datum in $U$ that contacts $\cald$, then there exists a datum $\zeta$ in $\d$ such that $\delta =k(\zeta)$.
\end{enumerate}
\end{prop}
In honor of Lempert's seminal theorem \cite{lem81} we adopt the following definition.
\begin{defin}\label{predef50} 
A domain $U$ in $\c^n$ is a \emph{Lempert domain}
\index{Lempert!domain}
 if 
\begin{enumerate}[\rm (1)]
\item $U$ is weakly hyperbolic
\item  $U$ is taut and 
\item $\car{\delta}=\kob{\delta}$ for every nondegenerate datum $\delta$ in $U$.
\end{enumerate}
\end{defin}
 Equivalently, $U$ is a Lempert domain if  $U$ is hyperbolic in the sense of Kobayashi \cite{kob98}, meaning that the topology induced by the Kobayashi pseudodistance on $U$ is the Euclidean topology, $U$ is taut and every datum in $U$  contacts  a complex geodesic.
\begin{lem}\label{prelem20}
Let $\Omega$ be a Lempert domain in $\c^n$ and $S$ be a  subset of $\Omega$ having nonempty interior.   Suppose that, for all $\mu_1,\mu_2 \in S$ and for every complex geodesic $\mathcal{D}$ in $\Omega$ containing $\mu_1$ and $\mu_2$, we have $\mathcal{D} \subseteq S$. Then $S=\Omega$.
\end{lem}
\begin{proof}
Fix $\mu \in \Omega$. Choose an interior point $\mu_1$ of $S$ and let $\mathcal{D}_1$ be the complex geodesic passing through $\mu$ and $\mu_1$. Since  $\mu_1$ is an interior point of $S$, there exists a point $\mu_2 \in S\cap\cald_1$ such that $\mu_2 \neq  \mu_1$.  Thus $\mu_1, \mu_2 \in S\cap\cald_1$,  and so $\mathcal{D}_1$ is a geodesic passing through $\mu_1$ and $\mu_2$. It follows by assumption that $\mathcal{D}_1 \subseteq S$. In particular $\mu \in S$.
\end{proof}

%% file: charactbidisc.tex
\section{Characterizations of the disc and bidisc}\label{charactbidisc}

The following is the simplest result on the characterization of a domain through its Carath\'eodory extremal functions.
\begin{theorem}\label{chardisc}
Let $\Omega$ be a Lempert domain.  If there is a function $\Phi\in\d(\Omega)$ such that $\{\Phi\}$ is a universal set for the Carath\'eodory problem on $\Omega$ then $\Omega$ is isomorphic to $\d$.
\end{theorem}
For the proof we shall invoke the following observation.
\begin{lemma}
Let $\Omega$ be a Lempert domain.  If $\Phi\in\d(\Omega)$ and $\Phi$ solves $\Car\de$ for some \nd datum $\de$ in $\Omega$ then $\Phi$ is surjective.
\end{lemma}
\begin{proof}
By the tautness of $\Omega$, there exists a function $g\in\Omega(\d)$ that solves $\Kob \de$.  Hence there is a datum $\zeta$ in $\d$ such that $g(\zeta)=\de$ and
\[
|\zeta|=\kob{\de}^\Omega=\car{\de}^\Omega=|\Phi(\de)|= |\Phi\circ g(\zeta)|.
\]
Thus $\Phi\circ g$ is a holomorphic self-map of $\d$ that preserves the modulus of a \nd datum in $\d$.  Hence  $\Phi\circ g$ is an automorphism of $\d$, from which it follows that $\Phi$ is surjective.
\end{proof}
\begin{proof}[Proof of Theorem \ref{chardisc}]
Consider any \nd discrete datum $\de$ in $\Omega$.  By the definition of Lempert domains, $\Omega$ is weakly hyperbolic and so $\car{\de} > 0$.  Since $\Phi$ solves $\Car \de$,
\[
|\Phi(\de)|=\car{\de} >0.
\]
It follows that $\Phi:\Omega\to\d$ is injective.

$\Phi$ is also surjective, by the preceding lemma.  Hence $\Phi$ is a holomorphic bijection, and therefore an isomorphism between $\Omega$ and $\d$.
\end{proof}

More remarkably, the bidisc can be characterized up to isomorphism among Lempert domains by the existence of a universal set for the Carath\'eodory problem comprising two functions.

For the proof of this statement
we need some ideas from \cite[Section 1]{agmc_vn}.  We shall say that a discrete datum $\la=(z,w)$ in $\d^2$ is {\em balanced} if it is \nd and
\index{balanced}
\[
|(z^1,w^1)|= |(z^2,w^2)|.
\]
When this equation holds there is a unique automorphism $m$ of $\d$ such that $m(z^1)=z^2$ and $m(w^1)=w^2$.
Furthermore there is a unique complex geodesic $D_\la$ in $\d^2$ that is contacted by $\la$, to wit
\[
D_\la=\{(\zeta,m(\zeta)):\zeta\in\d\}.
\]

A subset $V$ of $\d^2$ is said to be {\em balanced} if, for every balanced datum $\la$ in $\d^2$ that contacts $V$, $D_\la \subseteq V$.
According to \cite[Proposition 4.10]{agmc_vn}, if $B$ is a balanced subset of $\d^2$ that is contacted by a balanced datum $\la$, then either $B=D_\la$ or $B=\d^2$.

\begin{thm}\label{thm10}
Let $\Omega$ be a Lempert domain.  $\Omega$ is biholomorphic to $\d^2$ if and only if there is a minimal universal set for the Carath\'eodory problem on $\Omega$ consisting of two functions.
\end{thm}

\begin{proof} Firstly let us show that
the set of the two co-ordinate functions is a minimal universal set for the Carath\'eodory problem on $\d^2$. Let $F^j(\lambda)= \lambda^j$ for $ \lambda= (\lambda^1,\lambda^2) \in \d^2$. 

Consider any \nd discrete datum $(\la_1,\la_2)$ in $\d^2$, where 
 $\lambda_1= (\lambda_1^1,\lambda_1^2)$ and $\lambda_2= (\lambda_2^1,\lambda_2^2)$. 
We have
\[
\car{(\lambda_1, \lambda_2)}=
\sup_{F\in \d(\d^2)} |(F(\lambda_1), F(\lambda_2))| \ge 
\max \{|(\lambda_1^1,\lambda_2^1)|, |(\lambda_1^2,\lambda_2^2)| \}.
\]
Suppose that
\beq\label{type1}
|(\lambda_1^1,\lambda_2^1)|\ge |(\lambda_1^2,\lambda_2^2)|,
\eeq
so that 
\[
\car{(\lambda_1, \lambda_2)}=|(\lambda_1^1,\lambda_2^1)|.
\]

By inequality \eqref{type1} and the Schwarz-Pick Lemma, there exists 
 $f \in \d(\d)$ such that 
\beq\label{9.1}
f(\lambda_1^1)=\lambda^2_1\quad\mbox{ and }\quad f(\lambda_2^1)=\lambda^2_2.
\eeq
Consider any extremal function $F \in \d(\d^2)$ for $\Car (\la_1,\la_2)$.
Then
\[
 \car{(\lambda_1, \lambda_2)}=|(F(\lambda_1),F (\lambda_2))|.
\]

Define $\tilde{F} \in \d(\d^2)$  by $\tilde F(z)=F(z,f(z))$ for $z\in\d$.  Then, for $i=1,2$,
\[
\tilde F(\la_i^1)=F(\la_i^1,f(\la_i^1))=F(\la_i^1,\la_i^2)=F(\la_i).
\]
Again by the Schwarz-Pick Lemma,
\[
\left|\left(\tilde F(\la^1_1), \tilde F(\la^1_2)\right)\right| \leq |(\la^1_1,\la^1_2)|.
\]
Hence
\begin{align*}
|(\la_1,\la_2)|^{\d^2}_{\mathrm {car}} &\leq \max \{ |(\la^1_1,\la^1_2)|,  |(\la^2_1,\la^2_2)|\}\\
	&=\max\{ |F^1((\la_1,\la_2))|, |F^2((\la_1,\la_2))|\},
\end{align*}
so that either $F^1$ or $F^2$ is extremal for $\Car (\la_1,\la_2)$.

A similar proof applies to infinitesimal datums in $\d^2$.
Thus $\{F^1,F^2\}$ is a universal set for the Carath\'eodory problem on $\d^2$.

Clearly, neither $\{F^1\}$ nor $\{F^2\}$ is a universal set for the Carath\'eodory problem on $\d^2$.
Hence $\{F^1,F^2\}$ is a minimal universal set for the Carath\'eodory problem on $\d^2$.

To prove the converse, let $\{\Psi_1,\Psi_2\}$ be a universal set for the Carath\'eodory problem on $\Omega$ 
and define a holomorphic mapping $\Psi$ on $\Omega$ by the formula
\[
\Psi(\mu) = (\Psi_1(\mu), \Psi_2(\mu)).
\]
Since $\Psi_1,\Psi_2 \in \d(\Omega)$, clearly $\Psi$ maps $\Omega$ to $\d^2$.
That $\Omega$ is biholomorphic to $\d^2$ will follow if it is shown that $\Psi$ is bijective from $\Omega$ to $\d^2$,
since every bijective holomorphic map between domains has a holomorphic inverse (for example, \cite[Chapter 10, Exercise 37]{kr}).

To see that $\Psi$ is injective, consider distinct points  $\mu_1$ and $\mu_2$ in $\Omega$. As $\Omega$ is a Lempert domain, $\Omega$ is weakly hyperbolic, and so
\begin{align*}
0 &< \kob{(\mu_1,\mu_2)}=\car{(\mu_1,\mu_2)} \\
	&=\max \left\{|(\Psi_1(\mu_1),\Psi_1(\mu_2))|,
 |(\Psi_2(\mu_1),\Psi_2(\mu_2))|\right\}.
\end{align*}
Hence either $\Psi_1(\mu_1)\not= \Psi_1(\mu_2)$ or $\Psi_2(\mu_1)\not= \Psi_2(\mu_2)$, so that $\Psi(\mu_1) \not= \Psi(\mu_2)$.

To prove that $\Psi$ is surjective, 
we first show that there is a \nd datum $(\mu,\nu)$ in $\Omega$ such that $\Psi((\mu,\nu))$ is a balanced datum in $\d^2$.

Since $\{\Psi_1,\Psi_2\}$ is a {\em minimal} universal set for $\Omega$, $\{\Psi_2\}$ is not a universal set.  Therefore there is a \nd datum $(\mu_1,\nu_1)$ in $\Omega$ for which $\Psi_1$ is a Carath\'eodory extremal but $\Psi_2$ is not.   Thus
\[
\left|\left(\Psi_1(\mu_1),\Psi_1(\nu_1)\right)\right|\quad  > \quad \left|\left(\Psi_2(\mu_1),\Psi_2(\nu_1)\right) \right|.
\]
Similarly, there is a \nd datum $(\mu_2,\nu_2)$ in $\Omega$ such that
\[
\left|\left(\Psi_1(\mu_2),\Psi_1(\nu_2)\right)\right|\quad  < \quad \left|\left(\Psi_2(\mu_2),\Psi_2(\nu_2)\right) \right|.
\]

As $\Omega$ is a Lempert domain, every pair of points in $\Omega$ lies in an analytic disc in $\Omega$, which implies that $\Omega$ is connected.  Thus $\Omega\times\Omega$ is connected, and consequently the set of \nd~ datums in $\Omega$,
\[
\mathrm{ndd}(\Omega) \df (\Omega\times\Omega) \setminus \{(\mu,\mu):\mu\in\Omega\},
\]
is connected (note that the diagonal set $\{(\mu,\mu)\}$ lies in a subspace of real codimension at least two in $\Omega\times\Omega$).
Consider a continuous path
\[
\ga = (\ga_1,\ga_2): [0,1] \to \mathrm{ndd}(\Omega)
\]
such that $\ga(0)=(\mu_1,\nu_1)$ and $\ga(1)=(\mu_2,\nu_2)$.  Then
\[
f(t) \df \left|\left( \Psi_1\circ \ga_1(t), \Psi_1\circ \ga_2(t)\right)\right| - \left|\left(\Psi_2\circ\ga_1(t),\Psi_2\circ\ga_2(t)\right)\right|
\]
depends continuously on $t$ for $0\leq t\leq 1$, is strictly positive at $t=0$ and strictly negative at $t=1$.  Hence there exists $t_0\in(0,1)$ such that $f(t_0)=0$.   Then $(\mu_0,\nu_0)\df \ga(t_0)$ is a \nd datum in $\Omega$ and
\be\label{abaldat}
\Psi((\mu_0,\nu_0)) \mbox{ is a balanced datum in } \d^2.
\ee

Now we show that $\Psi(\Omega)$ is a balanced set in $\d^2$.   Consider any balanced discrete datum $\la=(z,w)$ in $\d^2$ that contacts $\Psi(\Omega)$.   We wish to prove that $D_\lambda$, the balanced disc in $\d^2$ passing through $z$ and $w$, is contained in $ \Psi(\Omega)$.

Since $z,w\in\Psi(\Omega)$, we may choose points $\mu,\nu\in\Omega$ such that $z=\Psi(\mu)$ and $w=\Psi(\nu)$. 
Since $\{\Psi_1,\Psi_2\}$ is universal for the Carath\'eodory problem on $\Omega$,
\begin{align*}
\car{(\mu,\nu)}^\Omega &= \max\left\{ |(\Psi_1(\mu),\Psi_1(\nu))|,|(\Psi_2(\mu),\Psi_2(\nu))|\right\}\\
	&=\car{(\Psi(\mu),\Psi(\nu))}^{\d^2}.
\end{align*}

Since $\Omega$ is taut there exists a function $g\in\Omega(\d)$ that solves $\Kob((\mu,\nu))$.  That is, there exist $\al,\beta\in\d$ such that $|(\al,\beta)|=\kob{(\mu,\nu)}$ and $g(\al)=\mu, \, g(\beta)=\nu$.
We have now $\Psi\circ g \in \d^2(\d)$ while $\Psi\circ g(\al)=z, \, \Psi\circ g(\beta)=w$ and
\begin{align*}
|(\al,\beta)|&=\kob{(\mu,\nu)}^\Omega=\car{(\mu,\nu)}^\Omega=\car{(\Psi(\mu),\Psi(\nu)}^{\d^2}\\
	&=\kob{(\Psi\circ g(\al),\Psi\circ g(\beta))}^{\d^2}.
\end{align*}
Thus $\Psi\circ g$ solves $\Kob(\la)$.  Since $\la$ is a balanced datum in $\d^2$, there is a {\em unique} complex geodesic in $\d^2$ contacted by $\la$, and therefore 
\[
\Psi\circ g (\d)=D_\la.
\]
Hence $D_\la \subseteq \Psi(\Omega)$.  It follows that $\Psi(\Omega)$ is a balanced set in $\d^2$.

Let 
\begin{align*}
\la_0 &= \Psi((\mu_0,\nu_0)) \\
	&= (\Psi_1,\Psi_2)((\mu_0,\nu_0))\\
	&= \left( ( \Psi_1(\mu_0),\Psi_2(\mu_0)), (\Psi_1(\nu_0),\Psi_2(\nu_0))\right).
\end{align*}
Since, by equation \eqref{abaldat}, $\Psi(\la_0)$ is a balanced datum that contacts $\Psi(\Omega)$,
 \cite[Proposition 4.10]{agmc_vn} applies and yields the conclusion that either $\Psi(\Omega)= D_{\la_0}$ or $\Psi(\Omega)= \d^2$.

Suppose that 
\[
\Psi(\Omega)= D_{\la_0}=\{(\zeta,m(\zeta)):\zeta\in\d\}.
\]
Then $\Psi_2=m\circ\Psi_1$.  Consequently, $\{\Psi_1\}$ is also a universal set for the Carath\'eodory problem on $\Omega$, contrary to the minimality of $\{\Psi_1,\Psi_2\}$.   Hence
$\Psi(\Omega)\neq D_{\la_0}$, and therefore $\Psi(\Omega)= \d^2$.

Thus $\Psi:\Omega\to\d^2$ is a bijective holomorphic map between domains, and is therefore an isomorphism.
\end{proof}

An interesting feature of the two theorems in this section is that the dimensions of the domains are obtained as consequences of assumptions about universal sets for the Carath\'eodory problems on the domains.

%% file: charGcara.tex
\section{A characterization of $G$ via Carath\'eodory extremals} \label{charactcara}

\subsection{Extremal problems and geodesics in $G$} 
In this subsection we recall some known facts about complex geodesics in the symmetrized bidisc.

\begin{thm}\label{prethm10}
For every $\omega \in \t$ define the rational function $\Phi_\omega$ by
\be\label{pre20}
\Phi_\omega(s^1,s^2) = \frac{2\omega s^2 -s^1}{2-\omega s^1}.
\ee
The set $\;\mathcal{C} = \{\Phi_\omega: \omega \in \t\}$ is a minimal universal set for the Carath\'eodory problem on $G$.
\end{thm}
The fact that $\mathcal C$ is a universal set is proved in \cite[Theorem 1.1 and Corollary 4.3]{ay2004}.
The minimality of $\mathcal{C}$ follows from the following fact.
\begin{lem}\label{unomega}
For every point $\tau\in\t$ there exists a \nd datum $\de$ in $G$ such that, for $\omega \in\t$, $\Phi_\omega$ solves $\Car\de$ if and only if $\omega=\tau$.
\end{lem}
\begin{proof}
This fact is shown in \cite{aly2016}.  Here is one way to construct such a $\delta$ \cite[Section 1]{aly2018b}.  Choose an automorphism $m$ of $\d$ having $\tau$ as its unique fixed point in the closed unit disc.  Let
\[
h(z) = (z+m(z),zm(z))
\]
for all $z\in\d$, and let $\de$ be the infinitesimal datum $(h(z_0),h'(z_0))$ for some $z_0\in\d$.  Then $\de$ has the required property.
\end{proof}
For $\tau$ and $\de$ as in the lemma, there is no solution of $\Car\de$ in $\mathcal{C} \setminus \{\Phi_\tau\}$.
Therefore, no proper subset of $\mathcal{C}$ is universal for the Carath\'eodory problem on $G$.

The following uniqueness result for the Kobayashi problem in $G$ is proved in \cite[Theorem 0.3]{AY06} for discrete datums and in \cite[Theorem A.10]{aly2016} for infinitesimal datums.
\begin{thm}\label{prethm20}
If $\delta$ is a nondegenerate datum in $G$, then the solution to $\Kob\delta$ is essentially unique.
\end{thm}
A surprising fact about $G$ is that Lempert's conclusion remains true despite the fact that $G$ is not convex (nor even biholomorphic to a convex domain  \cite{cos04}). To be specific, the following result is true 
 (\cite[Corollary 5.7]{ay2004} in the discrete case and \cite[Proposition 11.1.7]{jp} in the infinitesimal case).

\begin{thm}\label{prethm30}
$G$ is a Lempert domain, that is, if $\delta$ is a nondegenerate datum in $G$, then there exists a complex geodesic $\cald$ such that $\delta$ contacts $\cald$.
\end{thm}

On combining the last two theorems we deduce: 
\begin{corollary}\label{prethm40}
For every \nd datum $\de$ in $G$ there is a unique complex geodesic in $G$ that is contacted by $\de$.
\end{corollary}
As a consequence of these theorems we may unambiguously attach to each \nd datum in $G$ a unique complex geodesic.
\begin{defin}\label{predef60}
For any \nd datum $\delta$  in $G$,  $\cald_\delta$
\index{$\cald_\delta$}
 denotes the unique complex geodesic in $G$ that is contacted by $\delta$.
\end{defin}
\subsection{A characterization of $G$}
In this subsection we characterize $G$ in terms of the possession of a universal set for the Carath\'eodory problem of the same algebraic form as the universal set for $G$ described in Theorem \ref{prethm10}.

\begin{defin}\label{cardef10}
A domain $\Omega$ in $\c^2$ {\em has a $G$-like universal set} if there exist $s,p \in \c(\Omega)$ such that $\{\Psi_\omega\}_{\omega \in \t}$ is a universal set for the Carath\'eodory extremal problem on $\Omega$, where for each $\omega \in \t$, $\Psi_\omega$ is defined by
\[
\Psi_\omega(\mu)=\frac{2\omega p(\mu)-s(\mu)}{2-\omega s(\mu)}  \qquad \mbox{for all } \mu \in \Omega.
\]
\end{defin}

\begin{thm}\label{carthm10}
 If $\Omega$ is a domain in $\c^2$, then $\Omega$ is biholomorphically equivalent to the symmetrized bidisc if and only if $\Omega$ is a Lempert domain and $\Omega$ has a $G$-like universal set.
\end{thm}

\begin{proof} Clearly, if $\Omega$ is a domain in $\c^2$ and $F:\Omega \to G$ is biholomorphic, then, as $G$ is a Lempert domain, so also is $\Omega$. Furthermore, as $\set{\Phi_\omega}_{\omega \in \t}$ is a universal set for the Carath\'eodory extremal problem on $G$, if we define $\Psi_\omega=\Phi_\omega \circ F$ for all $\omega \in \t$, then $\set{\Psi_\omega}_{\omega \in \t}$ is a universal set for the Carath\'eodory extremal problem on $\Omega$. This proves that $\Omega$ has a $G$-like universal set.

Now assume that $\Omega$ is a Lempert domain in $\c^2$ and has a $G$-like universal set as in Definition \ref{cardef10}. Let $F=(s,p)$, so that $F$ is a holomorphic map $\Omega \to \c^2$.  We shall show that $F$ is a biholomorphic mapping of $\Omega$ onto $G$.

To see that $F(\Omega)\subset G$, consider $\mu\in\Omega$. Since $\Psi_\omega$ maps  $\Omega$ to $\d$,
\[
\left| \frac{2\omega p(\mu)-s(\mu)}{2-\omega s(\mu)}\right| < 1
\]
 for all $\omega \in \t$.
Hence
\[
|2\omega p(\mu)-s(\mu)|^2< |2-\omega s(\mu)|^2.
\]
Expand this relation to obtain
\[
 \re\left( \omega(s(\mu) -\overline{s(\mu)}p(\mu))\right) <1-|p(\mu)|^2
\]
for all $\omega\in\t$.
Consequently,
\[
  |s(\mu) -\overline{s(\mu)}p(\mu)|<1-|p(\mu)|^2.
\]
This inequality is equivalent to the statement that $(s(\mu),p(\mu))\in G$ (for example \cite[Theorem 2.1]{ay2004}), so that  $F(\mu) \in G$ for all $\mu \in \Omega$. This proves that $F\in G(\Omega)$.

We next prove that $F$ is injective and unramified (that is, $F'(\mu)$ is nonsingular for all $\mu \in \Omega$). Fix a nondegenerate datum $\delta$ in $\Omega$. Since $\Psi_\omega = \Phi_\omega \circ F$, 
\[
\Psi_\omega(\delta) = \Phi_\omega(F(\delta))
\]
for each $\omega \in \t$. By Lemma \ref{prelem10}, there exists $\omega\in\t$ such that $\Psi_\omega(\de)$ is nondegenerate, that is, $\Phi_\omega(F(\de))$ is nondegenerate.  This fact in turn implies that $F(\delta)$ is nondegenerate. To summarize, we have proved that if $\delta$ is nondegenerate then $F(\delta)$ is nondegenerate. The case when $\delta$ is discrete yields the conclusion that $F$ is injective, and the case when   $\delta$ is infinitesimal implies that $F$ is unramified.

It remains to prove that $F$ is surjective. Note that if $\mu_1,\mu_2 \in \Omega$, then
\begin{align*}
\car{(\mu_1,\mu_2)}^\Omega &= \sup \set{\left|\big(\Psi_\omega(\mu_1),\Psi_\omega(\mu_2)\big)\right|}{\omega \in \t}\\
&=\sup \set{\left|\big(\Phi_\omega(F(\mu_1)),\Phi_\omega (F(\mu_2))\big)\right|}{\omega \in \t}\\
&=\car{(F(\mu_1),F(\mu_2))}^G.
\end{align*}
Thus $F:\Omega \to G$ is an isometry when $\Omega$ and $G$ are equipped with their respective Carath\'eodory (or Kobayashi) distances.

Consider a \nd discrete datum $\mu=(\mu_1,\mu_2)$ in $\Omega$ and let  $\la=(\la_1,\la_2)=F(\mu)=(F(\mu_1),F(\mu_2))$.
Since $F$ is an isometry,
\be\label{isokob}
\kob{\la}^G = \kob{\mu}^\Omega.
\ee
By Corollary \ref{prethm40} there is a unique complex geodesic $D_\la$ in $G$ contacted by  $\lambda$.
We wish to prove that $D_\lambda$ is contained in $\ran F$.  Choose $g\in \Omega(\d)$ that solves $\Kob\mu$ and then $\alpha_1,\alpha_2 \in \d$ such that $g(\alpha_1)=\mu_1$, $g(\alpha_2)=\mu_2$, so that 
\be\label{almu}
|(\al_1,\al_2)|= \kob{\mu}^\Omega.
\ee

 Then $F \circ g \in  G(\d)$ and
\[
(F \circ g) (\alpha_i) = \lambda_i,\qquad i=1,2.
\]
On combining equations \eqref{isokob} and \eqref{almu} we obtain the statement
\begin{align*}
|(\al_1,\al_2)|= \kob{\la}^G.
\end{align*}
Thus $F \circ g$ solves the Kobayashi extremal problem for $\la$ in $G$.  Now the complex geodesic $D_\la$ is the range
of any solution of the Kobayashi problem for the \nd datum $\la$ in $G$ (see for example \cite[Theorem 4.6]{aly2016}) and so
\[
D_\lambda = \ran (F \circ g) \subseteq \ran F,
\]
as was to be proved.

The fact that $F$ is unramified guarantees that $\ran F$ contains a nonempty open set in $G$.  By  Lemma \ref{prelem20}, $\ran F = G$. 
We have shown that $F:\Omega\to G$ is a bijective holomorphic map.
\end{proof}

%% file: uniqueness.tex
\section{Uniqueness of minimal universal sets}\label{uniqueness}
It is natural to ask whether, for a general Lempert domain $\Omega$, there is a {\em unique}
minimal universal set for the Carath\'eodory extremal problem on $\Omega$, up to an obvious
notion of equivalence: if $\calc$ is a universal set for $\Omega$, then so is
\[
\calc'\df\{m_\ph\circ\ph : \ph\in\calc\}
\]
where $m_\ph$ is an automorphism of $\d$ for every $\ph\in\calc$.  We regard $\calc$ and $\calc'$
as equivalent universal sets.

We do not know whether uniqueness (up to equivalence) holds for a general Lempert domain, but for the three domains
studied in this paper, it does hold.  

In the case that $\Omega=\d$, as we observed in Example \ref{ex1.2}, it follows from the Schwarz-Pick Lemma that any singleton set containing an automorphism of $\d$
is a universal set for the Carath\'eodory problem on $\d$.  Such a set is clearly minimal.     Conversely, let $\calc$ be a universal set for $\d$.
Then $\calc$ contains a Carath\'eodory extremal function $\ph$ for the discrete datum $(0, \half)$.  Again the Schwarz Lemma implies that $\ph$
is an automorphism $m$ of $\d$, and so, by minimality, $\calc=\{m\}$.

Consider the case that $\Omega=\d^2$.  By,  for example, \cite[page 293] {agmc_vn}, a universal set for the Carth\'eodory problem on $\d^2$ is 
$\{F_1,F_2\}$, where $F_j$ is the $j$th co-ordinate function, given by $F_j(z)=z^j$ for $j=1,2$.  This set is easily seen to be minimal.
Conversely, let $\calc$ be a minimal universal set for $\d^2$.  Then $\calc$ contains a Carath\'eodory extremal for the unbalanced discrete datum $( (0,0), (\half,0))$.  Now such a datum has a {\em unique} Carath\'eodory extremal function, up to equivalence, to wit, the co-ordinate function $F_1(z)= z^1$ (see, for example, \cite[Theorem 12.2] {jimjohnbook}).  Hence $\calc$ contains $m_1\circ F_1$ for some automorphism $m_1$ of $\d$.  Likewise, consideration of the unbalanced datum $((0,0),(0,\half))$ leads to the conclusion that $\calc$ contains $m_2\circ F_2$, where $F_2$ is the second co-ordinate function and $m_2$ is an automorphism of $\d$.  Hence
\[
\{m_1\circ F_1, m_2\circ F_2 \} \subseteq \calc.
\]
Since the left hand side of this inclusion is a universal set for $\d^2$, it follows by minimality that the inclusion holds with equality.  Thus $\{F_1,F_2\}$ is the unique universal set for $\d^2$, modulo equivalence.

By Theorem \ref{prethm10}, the set $\{\Phi_\omega:\omega\in\t\}$ is a minimal universal set for the Carath\'eodory extremal problem on $G$.
Consider any other minimal universal set $\calc$ for $G$.  Let $\tau\in\t$. By Lemma \ref{unomega} there is a \nd datum $\delta$ in $G$ such that $\Phi_\omega$ is a Carath\'eodory extremal for $\delta$ if and only if $\omega=\tau$.    By \cite[Corollary 2.8]{aly2018b}, every Carath\'eodory extremal function for $\de$ is of the form $\gamma\circ \Phi_{\tau}$ for some automorphism $\gamma$ of $\d$.
Hence, up to equivalence, $\calc$ contains $\Phi_{\tau}$, and since $\tau\in\t$ was arbitrary,
\[
\{\Phi_\omega:\omega\in\t\} \subseteq \calc.
\]
By minimality, $\calc$ is equivalent to $\{\Phi_\omega:\omega\in\t\}$.

%% file: concluding.tex
\section{Concluding remarks}
Our results show that a minimal universal set for the  Carath\'eodory extremal problem on a domain $\Omega$
provides significant information about $\Omega$.  It would also be of interest to describe
{\em all}  Carath\'eodory extremal functions on $\Omega$.
 As we mentioned in the introduction, even for the bidisc, the set of all
 Carath\'eodory extremal functions is large.  
This fact is pointed out by \L. Kosi\'nski and W. Zwonek in \cite{kos}.  They discuss the Carath\'eodory extremal functions for bounded convex domains, strongly linearly convex domains, the symmetrized bidisc and the tetrablock. They describe cases in which the Carath\'eodory extremal function for a particular pair of points is unique modulo automorphisms of $\d$ and analyse the relationship between this property and the {\em non-uniqueness} of complex geodesics through the points.  The authors observe that their results give an understanding of the phenomenon of the uniqueness of Carath\'eodory extremal functions, but state that in the non-unique case the form of the extremal functions is not well understood.
We have studied the nature of Carath\'eodory extremals in the case of the symmetrized bidisc \cite{aly2018b}.  We were able to describe large classes of Carath\'eodory extremal functions for datums of various types in $G$, both when they are unique and when they are non-unique. 

We have also found other characterizations of the symmetrized bidisc.
In \cite{ALY2018} we  characterize $G$ in terms of the action of the automorphism group of $G$.